\documentclass[12pt]{amsart}

\usepackage{amsmath,amssymb,amsfonts,mathrsfs,mathtools}
\usepackage{graphicx,xcolor}
\mathtoolsset{showonlyrefs}
\allowdisplaybreaks[1]

\usepackage{wasysym,stmaryrd}

\definecolor{darkgreen}{rgb}{0,0.5,0}
\definecolor{darkred}{rgb}{0.7,0,0}
\usepackage[colorlinks, 
citecolor=darkgreen, linkcolor=darkred]{hyperref}
\usepackage{soul,esint}
\usepackage{mathtools}
\mathtoolsset{showonlyrefs}
\newtheorem{thm}{Theorem}[section]
\newtheorem{lma}[thm]{Lemma}
\newtheorem{prop}[thm]{Proposition}
\newtheorem{cor}[thm]{Corollary}
\theoremstyle{definition}
\newtheorem{defn}[thm]{Definition}
\theoremstyle{remark}
\newtheorem{rem}[thm]{Remark}
\numberwithin{equation}{section}
\usepackage{latexsym, amssymb, amsmath}

\newcommand{\be}{\begin{equation}}
\newcommand{\ee}{\end{equation}}
\newcommand{\Ric}{\operatorname{Ric}}
\newcommand{\Rm}{\operatorname{Rm}}
\newcommand{\Id}{\mathrm{Id}}

\newcommand{\scal}{\operatorname{scal}}
\newcommand{\inj}{\operatorname{inj}}
\newcommand{\loc}{\mathrm{loc}}
\newcommand{\R}{\mathbb{R}}
\newcommand{\C}{\mathbb{C}}
\newcommand{\e}{\varepsilon}
\def\a{\alpha}
\def\b{\beta}
\newcommand{\wt}{\widetilde}

\newcommand{\wn}{\widetilde\nabla}
\newcommand{\CPIC}{\mathrm{C}_{\mathrm{PIC1}}}

\begin{document}

\title[Weak stability of Ricci flow with almost WPIC1 curvature]
{Weak stability of Ricci flow with almost WPIC1 curvature}

\author{Pak-Yeung Chan}
\address[Pak-Yeung Chan]{Department of Mathematics, National Tsing Hua University, Hsin-Chu, Taiwan}
\email{pychan@math.nthu.edu.tw}

\author{Lucas Lavoyer}
\address[Lucas Lavoyer]{Mathematisches Institut, Universit\"at M\"unster, 48149 M\"unster, Germany}
\email{lucas.lavoyer@uni-muenster.de}

\author{Man-Chun Lee}
\address[Man-Chun Lee]{Department of Mathematics, The Chinese University of Hong Kong, Shatin, Hong Kong, China}
\email{mclee@math.cuhk.edu.hk}

\date{\today}

\begin{abstract}
We prove weak $L^\infty$ stability for complete, non-collapsed Ricci flows with almost WPIC1 curvature and scaling-invariant curvature bounds. Small $L^\infty$ perturbations of the initial metric admit Ricci-DeTurck flows which remain uniformly close to the background flow and depend Lipschitz continuously on the initial data. This extends an earlier stability theorem of the first and third authors with Lai. 
\end{abstract}

\maketitle

\section{Introduction}\label{sec:introduction}

On a  manifold $M$, the Ricci flow defined on an open interval $ M\times  I$ is a one parameter family of smooth metrics $g(t),t\in I$ such that it evolves along its Ricci direction: 
\begin{equation}
\partial_t g(t)=-2\Ric(g(t)).
\end{equation}
It was found to be powerful in regularizing metrics, especially in studying manifolds with positive or non-negative curvature. We refer readers to \cite{Brendle10} for an overview on its applications in sphere Theorem.

In this work, we are interested in the stability aspect of Ricci flow. The weak stability of Ricci flow concerns the dependence of a smooth flow on perturbations of its initial metric with only weak regularity. From the PDE perspective, the Ricci flow is a weakly parabolic system due to the invariance under diffeomorphisms. To fix the gauge, it was discovered by DeTurck \cite{DeTurck} that the Ricci flow is diffeomorphic to a strictly parabolic system which by now is called the Ricci-DeTurck flow. Precisely, given a smooth family of metrics  $\wt g(t)$ on a manifold $M$, the Ricci-DeTurck flow with respect to $\wt g(t)$ is a family of metrics $g(t)$ on $M$ satisfying
\begin{equation}\label{eqn:RDF0}
\left\{
\begin{aligned}
\partial_tg_{ij}&=-2\Ric(g)_{ij}+\nabla_iW_j+\nabla_jW_i,\\
W^k&=g^{ij}(\Gamma_{ij}^k(g)-\Gamma_{ij}^k(\wt g)).
\end{aligned}
\right.
\end{equation}

When $\tilde g(t)\equiv \tilde g$ is chosen to be a fixed metric with bounded geometry, the regularization of continuous metrics by the Ricci-DeTurck flow was first considered by Simon \cite{Simon2002}, who constructed short-time solutions for initial metrics sufficiently close to $\tilde g$ in $L^\infty$. In the case $\tilde g(t)\equiv g_{euc}$ on $\mathbb{R}^n$, Koch--Lamm \cite{KochLamm2012} established a global existence and uniqueness theory for sufficiently small $L^\infty$ perturbations of $g_{euc}$, together with analytic dependence on the initial data in suitable scaling-invariant spaces. It was later extended by  Burkhardt-Guim \cite{BurkhardtGuim2019} to compact manifolds; see also the recent work of Cai-Wang \cite{CaiWang2026}. There are also regularization results for initial metrics which need not be continuous, see the work of  Lamm-Simon \cite{LammSimon2023} for metrics in  $L^\infty\cap W^{2,2}$ in dimension four, and Chu and the third author \cite{ChuLee2025} for metrics in $L^\infty\cap W^{1,n}$. We also refer readers to \cite{Bamler2016,BurkhardtGuim2026,LitzingerSimon2025,ChuLeeWan,Lee2026} for the applications of stability in questions related to scalar curvature.

Motivated by question on smoothing singular spaces,  the first and third authors, together with Lai, started considering the stability of Ricci flow against a background Ricci flow $\tilde g(t)$ coming out of a cone $C(X)$ with non-negative curvature in a weak sense \cite{ChanLaiLee2025}. Those Ricci flows satisfy a bound of the form $|\Rm(\wt g(t))|\leq\alpha t^{-1},$ which does not give a uniform curvature bound down to time zero. Precisely, a weak stability theorem for backgrounds satisfying a curvature operator lower bound and scaling-invariant curvature and injectivity radius estimates is established in \cite[Section 4]{ChanLaiLee2025}. The purpose of this work is to extend this theorem to backgrounds with almost WPIC1 curvature. The situation where the background Ricci flow $\tilde g(t)$ has a stronger curvature asymptotic was studied earlier by Deruelle--Lamm \cite{DeruelleLamm2017}. 

We first recall the WPIC1 and WPIC2 curvature conditions, following \cite{ToppingPIC1} (see also \cite{Brendle10} for other equivalent definitions). Let $\mathcal{C}_B$ denote the space of algebraic curvature tensors on $\R^n$. We use the convention that $R(e_i,e_j,e_i,e_j)$ is the sectional curvature of the plane $\Sigma:=\mathrm{span}\{e_i,e_j\}$, and extend $R$ complex linearly to $\Lambda^2(\C^n)$. The Euclidean inner product also extends complex bilinearly to this space.

\begin{defn}\label{defn:PIC1}
An algebraic curvature tensor $R\in\mathcal{C}_B$ belongs to $\CPIC$ if
\begin{equation}
R(\omega,\bar\omega)\geq 0
\end{equation}
for every simple, isotropic bivector $\omega\in\Lambda^2(\C^n)\setminus\{0\}$, where isotropic means $\langle\omega,\omega\rangle=0$ for the complex bilinear inner product. Similarly, 
$R\in\mathcal{C}_B$ belongs to $\mathrm{C}_{\mathrm{PIC2}}$ if for any $a,b\in C^n$ with $a\wedge b\neq 0$,
\begin{equation}
R(a\wedge b,\bar a\wedge \bar b)\geq 0,
\end{equation}

A Riemannian metric has weakly PIC1 (weakly PIC2) curvature, abbreviated WPIC1 (resp. WPIC2), if $\Rm(g)\in \CPIC$ (resp. $\mathrm{C}_{\mathrm{PIC2}}$) at every point. Similarly we define strictly PIC1 and strictly PIC2 if $\Rm(g)$ lies inside the interior of the corresponding cone.
\end{defn}

We write $\Id:=\frac12 g\owedge g$ for the algebraic curvature tensor of constant sectional curvature one. The condition
\begin{equation}\label{eqn:almost-PIC1}
\Rm(g)+\kappa\cdot \Id\in\CPIC
\end{equation}
will be called almost WPIC1 curvature, with lower-bound parameter $\kappa\geq0$. In dimension three, $\CPIC$ is precisely the cone of nonnegative Ricci curvature, so \eqref{eqn:almost-PIC1} is equivalent to $\Ric(g)\geq-2\kappa g$.

We now state the main result of this paper.
\begin{thm}\label{thm:intro-main}
For every $n\geq3$ and $\a>0$, there exist $\bar\e_0(n,\a)>0$ and $\Lambda(n,\a)\geq1$ with the following property. Let $\wt g(t)$ be a smooth complete Ricci flow on $M\times[0,T]$, with $\wt g(0)=\wt g_0$, satisfying
\begin{equation} 
|\Rm(\wt g(t))|\leq\a t^{-1},\qquad
\inj(\wt g(t))\geq\sqrt{\a^{-1}t},\qquad
\Rm(\wt g(t))+\Id\in\CPIC.
\end{equation}
on $M\times(0,T]$. If $g_0$ is an $L^\infty$ metric such that
\begin{equation}
\|g_0-\wt g_0\|_{L^\infty(M,\wt g_0)}<\bar\e_0,
\end{equation}
then there is a smooth complete solution $g(t)$ to \eqref{eqn:RDF0} on $M\times(0,T\wedge1]$ with initial data $g_0$ and
\begin{enumerate}
\item[(a)]
$\displaystyle \|g(t)-\wt g(t)\|_{L^\infty(M,\wt g(t))}
\leq\Lambda\|g_0-\wt g_0\|_{L^\infty(M,\wt g_0)}$;
\item[(b)] if $g_0$ is continuous on an open set $\Omega\subset M$, then $g(t)\to g_0$ in $C^0_{\loc}(\Omega)$ as $t\downarrow0$.
\end{enumerate}
Moreover, if $g(t),\hat g(t)$ are two solutions with initial data $g_0,\hat g_0$ and
\begin{equation} 
\sup_{0<t\leq T\wedge1}\max\big\{
\|g(t)-\wt g(t)\|_{L^\infty},
\|\hat g(t)-\wt g(t)\|_{L^\infty}\big\}\leq\Lambda\bar\e_0,
\end{equation}
then
\begin{equation} 
\sup_{0<t\leq T\wedge1}\|g(t)-\hat g(t)\|_{L^\infty(M,\wt g(t))}
\leq\Lambda\|g_0-\hat g_0\|_{L^\infty(M,\wt g_0)}.
\end{equation}
In particular, the solution is unique in this class with its prescribed initial data.
\end{thm}

Here by {\sl initial data}, we mean $g(t)$ converges to $g_0$ in $L^p_{loc}$ sense  as $t\to 0^+$, for any $1\leq p<\infty$. In particular, it shows that the Ricci flows constructed in \cite{SimonToppingGT,Lai2019} are weakly stable. 

The core of our method is an improvement of controlling the linear curvature term appearing in the Lichnerowicz laplacian, based on the work of Zhao--Zhu \cite{ZhaoZhu22}. In the course of the proof, we obtain the strict linear stability and dynamic stability of expanding gradient Ricci soliton with WPIC1, extending a previous result by Deruelle--Lamm \cite{DeruelleLamm2017} under the setting of nonnegative sectional curvature, see Theorem~\ref{cor-strict-stab-PIC10}. 

One of the important examples of $\tilde g(t)$ in Theorem~\ref{thm:intro-main} is the smoothing of a metric cone $C(X)$ where $X=(\mathbb{S}^{n-1},h)$ is a smooth metric on $\mathbb{S}^{n-1}$. This models singular manifolds with isolated conical singularities. The Ricci flow smoothing of a metric cone $C(X)$ can be often shown to inherit the self-similar behaviour of the cone, and the solution is given by a so called expanding Ricci soliton. Recall that a triple $(M,g,f)$ consisting of a complete Riemannian manifold $(M,g)$ and a function $f\in C^\infty(M)$ is said to be an expanding gradient Ricci soliton (or an expander, for simplicity) if it satisfies
\[
\Ric(g)=\nabla^2f-\frac{1}{2}g.
\]

In \cite{Deruelle2016}, Deruelle first gave a PDE approach to the construction of expander smoothing starting from $C(X)$ when $X=(\mathbb{S}^{n-1},h)$ with $\Rm(h)\geq \mathrm{Id}$. Using Brendle's Harnack inequality \cite{Brendle09}, it is not hard to show that the same holds if we have only $\Rm(h)-\mathrm{Id}\in \mathrm{C}_{\mathrm{PIC2}}$. In the course of our proof, we show how our stability result yields the existence of an expander smoothing of such cones. Precisely, we show 
\begin{thm} 
   Let $(\mathbb{S}^{n-1},h)$ be a smooth closed manifold with $\Rm(h)-\Id_h\in \mathrm{C}_{\mathrm{PIC2}}$. Then there exists an expanding gradient Ricci soliton $(M,g,f,p)$ with $\Rm(g)\in \mathrm{C}_{\mathrm{PIC2}}$  asymptotic to the cone $C(\mathbb{S}^{n-1})$ over $(\mathbb{S}^{n-1},h)$, in the distance sense. Moreover if in addition $\Rm(h)-\Id_h\in \mathrm{int}\left(\mathrm{C}_{\mathrm{PIC2}}\right)$ somewhere on $\mathbb{S}^{n-1}$, then  $\Rm(g)\in\mathrm{int}(\mathrm{C}_{\mathrm{PIC2}})$.
\end{thm}

Furthermore, the solution is unique among solutions with the same metric cone at infinity in the distance sense, see Corollary~\ref{cor:uniq}.

Throughout the paper, $n\geq3$. Unless otherwise specified, tensor norms, inner products and covariant derivatives are taken with respect to $\wt g(t)$.  Constants denoted by $C_n$ depend only on $n$, and may change from line to line. We may replace $\a$ by $\max\{\a,1\}$ without further comment.

\subsection*{On the use of AI}
The writing and mathematical content of this work did not involve LLMs.

\subsection*{Acknowledgments}
The authors would like to thank Alix Deruelle, Felix Schulze, and Burkhard Wilking for useful discussion in the early stages of this work.  P.-Y. Chan is supported by the Yushan Young Fellow Program of the Ministry of Education (MOE), Taiwan (MOE-108-YSFMS-0004-012-P1), and by NSTC grants 113-2115-M-007-014-MY2 and 115-2628-M-007-006. L. Lavoyer is supported by the German Research Foundation – Project-ID 427320536 – SFB 1442, and by Germany’s Excellence Strategy EXC 2044/2 390685587, Mathematics Münster: Dynamics–Geometry–Structure. M.-C. Lee is supported by Hong Kong RGC grants No. 14300623 and No. 14304225, and an Asian Young Scientist Fellowship.

\section{An algebraic estimate under WPIC1}\label{sec:algebra}

In this section, we first collect some properties of the cone $\mathrm{C}_{\mathrm{PIC1}}$ of curvature operator with non-negative $1$-isotropic curvature. This will play a substantial role in establishing the weak stability of Ricci flow. The heart of the analysis relies on analyzing the Lichnerowicz Laplacian. In what follows, if $R$ is a curvature type operator, we use the convention that $R_{ijij}$ is the sectional curvature of
the plane spanned by $e_i,e_j$. For any symmetric $2$-tensor $h$, 
\begin{equation}\label{eqn:Rm-L_h}
    (\Rm(h))_{ij}:=\sum_{k,l}R_{ikjl}h_{kl}
\end{equation}
is the curvature term appearing in Lichnerowicz Laplacian. We also use $\Id=\frac12 g\owedge g$ to denote the curvature operator of standard sphere $\mathbb{S}^n$.

We will need the following algebraic inequality result obtained by Zhao--Zhu \cite{ZhaoZhu22}, which states that the curvature term \eqref{eqn:Rm-L_h} can be controlled by its scalar curvature, assuming $R\in\mathrm{C}_{\mathrm{PIC1}}$ and $n\geq 4$. We observe that the same holds for $n=3$ where $\mathrm{C}_{\mathrm{PIC1}}=\mathrm{C}_{\Ric}$.

\begin{prop}\cite{ZhaoZhu22}\label{prop-algebra under pic1}
    Let $n\geq 3$ and $R$ be an algebraic curvature tensor with $R\in \mathrm{C}_{\mathrm{PIC1}}$. Then, for any symmetric 2-tensor $h$, we have
    \begin{equation}
        \scal(R)|h|^2 - 2\langle \Rm(h), h\rangle \geq 0
    \end{equation}
where $\Rm(h)$ is defined in \eqref{eqn:Rm-L_h}.
\end{prop}
\begin{proof}
When $n\geq 4$, it follows directly from \cite[Proposition 1.6]{ZhaoZhu22}. We only need to prove $n=3$. Nevertheless we start with a general computation. We diagonalize $h$ with respect to $g$ so that $h(e_i, e_j) = \lambda_i \delta_{ij}$, and $\lambda_i \in \mathbb{R}$. If we denote the sectional curvature of $\Sigma=\mathrm{span}\{e_i,e_j \}$ to be $K_{ij}=R_{ijij}$, then  following the computation in \cite{DeruelleLamm2017} yields
\begin{equation}\label{eqn:ZhaoZhu}
    \begin{split}
 \scal(R)|h|^2 - 2\langle \Rm (h), h\rangle
  &=\left(\sum_{i<j} K_{ij}\right)\left(\sum_{k} \lambda_{k}^{2}\right) - 2\sum_{i<j} K_{ij} \lambda_{i} \lambda_{j}\\
  & = \sum_{i<j} K_{ij} \cdot \left(\sum_{k \neq i,j} \lambda_{k}^{2} + (\lambda_{i} - \lambda_{j})^{2}\right).
     \end{split}
\end{equation}

We now restrict ourself to $n=3$ so that 
\begin{equation}
    K_{ij}=\frac12 \left(\mu_i +\mu_j -\sum_{k\neq i,j}\mu_k\right)
\end{equation}
for all $i\neq j$, where $\mu_i:=\Ric_{ii}\geq 0$. Substituting it into \eqref{eqn:ZhaoZhu} yields 
\begin{equation}
    \begin{split}
 \scal(R)|h|^2 - 2\langle \Rm( h), h\rangle
  &=\frac12 \left(\mu_1+\mu_2-\mu_3 \right)\left( \lambda_3^2+(\lambda_1-\lambda_2)^2\right)\\
  &\quad +\frac12 \left(\mu_1+\mu_3-\mu_2 \right)\left( \lambda_2^2+(\lambda_1-\lambda_3)^2\right)\\
  &\quad +\frac12 \left(\mu_2+\mu_3-\mu_1 \right)\left( \lambda_1^2+(\lambda_3-\lambda_2)^2\right)\\
  &=\frac12 \mu_1 \left(  \lambda_2+\lambda_3-\lambda_1\right)^2+\frac12 \mu_2 \left(  \lambda_1+\lambda_3-\lambda_2\right)^2\\
  &\quad +\frac12 \mu_3 \left(  \lambda_2+\lambda_1-\lambda_3\right)^2\geq 0.
     \end{split}
\end{equation}
This completes the proof.
\end{proof}

As a consequence of Proposition~\ref{prop-algebra under pic1}, we use the stability result of Deruelle \cite{Der14b} to  extend   \cite[Lemma 5.15]{DeruelleLamm2017}, where nonnegative sectional curvature is assumed.

\begin{thm}\label{cor-strict-stab-PIC10}
Let $(M^n, g, \nabla^g f)$, $n \geq 3$, be an expanding gradient Ricci soliton with curvature operator in WPIC1. Then the bottom of the $L^2(d\mathrm{vol}_f)$-spectrum of the weighted Lichnerowicz operator is bounded from below by $n/2$, i.e.
\[
\lambda_1(-\Delta_f - 2 \Rm) \geq \frac{n}{2},
\]
where $\Delta_f:=\Delta+\nabla_{\nabla f} $ and $d\mathrm{vol}_{f}:=e^f d\mathrm{vol}_g$. In particular, $(M^n, g, \nabla^g f)$ is dynamically stable in the sense of \cite{Der14b}.
\end{thm}

\begin{proof}
    By the variational characterisation of $\lambda_1(-\Delta_f - 2\Rm)$, it suffices to prove that
\begin{align}
\inf_{h \in L^2(S^2T^*M, d\mu_f) \setminus \{0\}} \frac{\int_M \langle -\Delta_f h -2 \Rm(h), h \rangle \, d\mathrm{vol}_f}{\int_M |h|^2 d\mathrm{vol}_f} \geq \frac{n}{2}.
\end{align}

By \cite{MW2012, Der14b}, we know that for any $ \phi \in C_0^{\infty}(M)$, 
\begin{equation}
    \int_M \langle -\Delta_f \phi, \phi \rangle d\mathrm{vol}_f \geq \int_M \Delta f\cdot  \phi^2 d\mathrm{vol}_f,
\end{equation}
and $
\Delta f = \scal + n/2$.  For a symmetric 2-tensor with compact support $h$, the Kato inequality then yields
\begin{equation}\label{eqn:lamb}
\begin{split}
  &\quad   \int_M \langle -(\Delta_f + 2\Rm)h, h \rangle d\mathrm{vol}_f \\
  &\geq \int_M \langle (\scal - 2\Rm) h, h \rangle d\mathrm{vol}_f + \frac{n}{2} \int_M |h|^2 d\mathrm{vol}_f\geq \frac{n}{2} \int_M |h|^2 d\mathrm{vol}_f,
\end{split}
\end{equation}
where we have used Proposition \ref{prop-algebra under pic1}. From \cite{DengZhu15}, any gradient expander with nonnegative Ricci curvature must have bounded scalar curvature and proper potential $f$, so we can combine inequality \eqref{eqn:lamb}, \cite[Lemma A.2]{LeeTopping2025}, and \cite[Theorem 1.2]{Der14b} to show the dynamical stability of the expander.
\end{proof}

\section{A-priori estimates in space-time}\label{sec:RDF-esti}

Recall that, given a background Ricci flow $\wt g(t)$, the Ricci-DeTurck flow with respect to $\wt g(t)$ is a family of metrics $g(t)$ satisfying
\begin{equation}\label{eqn:RDF}
\left\{
\begin{aligned}
\partial_tg_{ij}&=-2\Ric(g)_{ij}+\nabla_iW_j+\nabla_jW_i,\\
W^k&=g^{ij}(\Gamma_{ij}^k-\wt\Gamma_{ij}^k).
\end{aligned}
\right.
\end{equation}
Above, $\nabla$ and $\Gamma$ are associated with $g(t)$, and the index of $W$ is lowered using $g(t)$. We use $\wn$ for the connection of $\wt g(t)$.

\subsection{The evolution inequality}

The algebraic estimates in Section~\ref{sec:algebra} allow us to extend  \cite[Lemma 4.3]{ChanLaiLee2025} to this weaker assumption.

\begin{lma}\label{lma:evo-h-L1}
Let $\tilde g(t)$ be a smooth family of metrics with $\Rm(\wt g)+\Id\in\CPIC$. For every $\gamma>0$, there exist $\e_0(\gamma,n),C_1(n)>0$ such that the following holds. Suppose $g=\wt g+h_1$ and $\hat g=\wt g+h_2$ are solutions to \eqref{eqn:RDF} so that 
\begin{equation}
|h_1|,|h_2|\leq \e_0,
\end{equation}
then the function $v:=(|h|^2+\sigma)^{1/2},$ where $0<\sigma<1$ and $h:=g-\hat g$, satisfies 
\begin{equation}\label{eqn:evo-general}
\begin{split}
&\quad \left(\frac{\partial}{\partial t}-\Delta_{\wt g}-\scal(\wt g)\right)v  \\
& \leq\gamma|\Rm(\wt g)|v+C_1|\Rm(\wt g)||h|^2
+C_1|h||\wn g|^2+C_1|\wn h|\big(|\wn g|+|\wn\hat g|\big)\\
&\quad+C_1v+\wn_k\big[(g^{kl}-\wt g^{kl})\wn_l v\big]
+\wn_p\big[v^{-1}\langle(g^{pq}-\hat g^{pq})\wn_q\hat g,h\rangle\big].
\end{split}
\end{equation}

Furthermore if $\hat g=\wt g$, there exists $\e_1(n)>0$ such that $|h|\leq\e_1$ implies
\begin{equation}\label{eqn:evo-single}
\begin{split}
&\quad \left(\frac{\partial}{\partial t}-\Delta_{\wt g}-\scal(\wt g)\right)v \\
&\leq C_1|\Rm(\wt g)||h|^2+C_1|\wn h|^2+C_1v
+\wn_k\big[(g^{kl}-\wt g^{kl})\wn_l v\big].
\end{split}
\end{equation}
\end{lma}
\begin{proof}

The proof is essentially the same as \cite[Lemma 4.3]{ChanLaiLee2025}, with the exception that we use Proposition~\ref{prop-algebra under pic1} to control the leading curvature term. By following the computation in \cite[Lemma 4.3]{ChanLaiLee2025}, we have 
\begin{equation}\label{eqn:h-evo}
\begin{split}
\partial_t |h|^2&=2 \langle \partial_t h,h\rangle \\
&\leq  2\langle g^{pq} \tilde \nabla_p\tilde\nabla_q g -\hat g^{pq} \tilde \nabla_p\tilde\nabla_q \hat g  ,h\rangle \\
&\quad +2\left(- g^{kl}  g_{ir}\tilde g^{rs} \tilde R_{jksl}- g^{kl}  g_{jr} \tilde g^{rs}\tilde R_{iksl}\right) h_{pq}\tilde g^{ip}\tilde g^{jq}\\
&\quad +2\left(\hat g^{kl} \hat  g_{ir}\tilde g^{rs} \tilde R_{jksl}+ \hat g^{kl}  \hat g_{jr} \tilde g^{rs}\tilde R_{iksl}\right) h_{pq}\tilde g^{ip}\tilde g^{jq}\\
&\quad +2h_{ik}\tilde R^k_j h_{pq}\tilde g^{ip}\tilde g^{jq}+2h_{jk}\tilde R^k_i h_{pq}\tilde g^{ip}\tilde g^{jq}\\
&\quad +C_n |h|^2 |\tilde \nabla g|^2+C_n |h||\tilde \nabla h|\left(|\tilde\nabla g|+|\tilde\nabla \hat g| \right).
\end{split}
\end{equation}
At a given point, choose a $\wt g$-orthonormal frame which diagonalizes $h$, with eigenvalues $\lambda_i$ so that  the terms involving Ricci curvature becomes 
\[2h_{ik}\tilde R^k_j h_{pq}\tilde g^{ip}\tilde g^{jq}+2h_{jk}\tilde R^k_i h_{pq}\tilde g^{ip}\tilde g^{jq}=4\langle\wt\Ric(h),h\rangle.
\]
The other curvature terms can be simplified as follows
\begin{eqnarray*}
&&-2\left(g^{kl}  g_{ir}\tilde g^{rs} \tilde R_{jksl}+g^{kl}  g_{jr} \tilde g^{rs}\tilde R_{iksl}\right) h_{pq}\tilde g^{ip}\tilde g^{jq}\\
&&+2\left(\hat g^{kl} \hat  g_{ir}\tilde g^{rs} \tilde R_{jksl}+ \hat g^{kl}  \hat g_{jr} \tilde g^{rs}\tilde R_{iksl}\right) h_{pq}\tilde g^{ip}\tilde g^{jq}\\
&=& 4\tilde R_{ikpl}\lambda_i(\hat g^{kl}\hat g_{ip}-g^{kl}g_{ip})\\
&=& -4\lambda_i^2 \tilde R_{ikil}\hat g^{kl}+4\lambda_i \lambda_p \tilde R_{ikrl}g^{kp}\hat g^{pl} g_{ir}\\
&=& -4\lambda_i^2 \tilde R_{ii}+4\lambda_i \lambda_k \tilde R_{ikik}+
h\ast h\ast \tilde R\ast \hat g^{-1}\ast h_2+h\ast h\ast\tilde R\ast g^{-1}\ast\hat g^{-1}\ast h_1\\
&\leq &-4\langle\wt\Ric(h),h\rangle+4\langle \wt\Rm(h),h\rangle + C_n |h|^2 (|h_2|+|h_1|)|\wt \Rm| 
\end{eqnarray*}

We now control the second curvature terms more carefully. Since $\wt\Rm+\mathrm{Id}\in \CPIC$, Proposition~\ref{prop-algebra under pic1} shows that  
\begin{equation}
\begin{split}
4\langle \wt\Rm(h),h\rangle&=4\langle (\wt\Rm+\mathrm{Id})(h),h\rangle-4|h|^2\\
&\leq 2\cdot \mathrm{scal}(\wt\Rm+\mathrm{Id}) \cdot |h|^2-4|h|^2\\
&=2\mathrm{scal}(\tilde g)\cdot |h|^2+2(n-2)(n+1)|h|^2.
\end{split}
\end{equation}

For the principal part, we have
\begin{equation}
\begin{split}
&\quad g^{pq}\wn_p\wn_qg-\hat g^{pq}\wn_p\wn_q\hat g \\
&=g^{pq}\wn_p\wn_qh
+\wn_p\big[(g^{pq}-\hat g^{pq})\wn_q\hat g\big]
-\wn_p(g^{pq}-\hat g^{pq})\cdot \wn_q\hat g.
\end{split}
\end{equation}
Together with $g^{-1}-\hat g^{-1}=-g^{-1}h\hat g^{-1}$, this gives the two divergence terms in \eqref{eqn:evo-general}. The remaining first-order terms can be controlled as in the proof of \cite[Lemma 4.3]{ChanLaiLee2025}.
\end{proof}

\subsection{A localized curvature integral}

An important second ingredient is a replacement of the spatial estimate of Petrunin \cite{Petrunin2008} by an estimate in space-time. 

\begin{lma}\label{lma:Petrunin}
Let $(M,\wt g(t))$ be a smooth Ricci flow, not necessarily complete, on $[0,T]$. Suppose that $x_0\in M$ and $r>0$ satisfy
\begin{enumerate}
\item[(i)] $B_{\wt g(0)}(x_0,4r)\Subset M$;
\item[(ii)] $|\Rm(\wt g(t))|\leq\a t^{-1}$ on $M\times(0,T]$;
\item[(iii)] $\Rm(\wt g(t))+r^{-2}\Id\in\CPIC$ on $M\times[0,T]$.
\end{enumerate}
There exist $\hat T_1(n,\a)\in(0,1]$ and $C_3(n)>0$ such that
\begin{equation}\label{eqn:spacetime-curvature}
\int_0^{(T\wedge\hat T_1)r^2}\int_{B_{\wt g(s)}(x_0,r)}
|\Rm(\wt g(s))|\,d\mathrm{vol}_{\wt g(s)}\,ds\leq C_3r^n.
\end{equation}
\end{lma}
\begin{proof}
We may assume $r=1$ by parabolic scaling. The cutoff construction in \cite[Lemma 7.1]{SimonToppingJDG} gives a function $\varphi$ with
\begin{equation}\label{eqn:evo-cutoff}
\begin{gathered}
0\leq\varphi\leq1,\qquad
(\partial_t-\Delta_{\wt g(t)})\varphi\leq0,\qquad
|\wn\varphi|^2\leq C_n\varphi,\\
\varphi\geq\tfrac12\ \text{on }B_{\wt g(t)}(x_0,1),\qquad
\operatorname{supp}\varphi(\cdot,t)\subset B_{\wt g(t)}(x_0,2),
\end{gathered}
\end{equation}
for $0\leq t\leq T\wedge\hat T_1$, for some $\hat T_1(n,\a)>0$.  The differential inequality holds in the barrier and hence in the distributional sense, see \cite[Appendix A]{MantegazzaMascellaniUraltsev}.

Consider the energy
\begin{equation}
E(t)=\int_M\varphi\,d\mathrm{vol}_{\wt g(t)},
\end{equation}
which Stoke Theorem implies
\begin{equation}
E'(t)\leq-\int_M\scal(\wt g(t))\varphi\,d\mathrm{vol}_{\wt g(t)}.
\end{equation}

Since $\wt\Rm+\mathrm{Id}\in \CPIC$, we have $\Ric(\tilde g)\geq -(n-1)$ so that 
\begin{equation}
E'(t)\leq n(n-1)E(t)
\end{equation}
for $t\in [0,T\wedge \hat T_1]$. Furthermore, volume comparison gives 
\begin{equation}\label{eqn:vol-bound}
E(t)\leq \mathrm{Vol}_{\tilde g(t)}\left(B_{\tilde g(t)}(x_0,2)\right)\leq C_n
\end{equation}
 for $t\in [0,T\wedge \hat T_1]$. Consequently,
\begin{equation}\label{eqn:weighted-scalar}
\begin{split}
\int_0^t\int_M \scal(\wt g(s))\varphi\,d\mathrm{vol}_{\wt g(s)}\,ds
&\leq E(0)-E(t)\leq C_n.
\end{split}
\end{equation}

From \cite[Lemma A.2]{LeeTopping2025}, there is $C_0(n)>0$ such that 
\begin{equation}\label{eqn:PIC1-contr-scal}
|\Rm(\tilde g)+\mathrm{Id}|\leq C_0 \left(\mathrm{scal}(\tilde g)+n(n-1) \right)
\end{equation}
By combining  \eqref{eqn:vol-bound}, \eqref{eqn:weighted-scalar} and \eqref{eqn:PIC1-contr-scal}, we obtain the space time integral estimate
\begin{equation}\label{eqn:weighted-curvature}
\int_0^t\int_M|\Rm(\wt g(s))|\varphi\,d\mathrm{vol}_{\wt g(s)}\,ds\leq C_n.
\end{equation}
Since $\varphi\geq1/2$ on $B_{\tilde g(t)}(x_0,1)$, this proves \eqref{eqn:spacetime-curvature}. 
\end{proof}

\subsection{Morrey estimates for metric perturbations}

Next, we prove a space time integral gradient estimate for the Ricci-DeTurck flow, generalizing \cite[Lemma 4.5]{ChanLaiLee2025} to the current setting.
\begin{lma}\label{lma:Morrey-grad}
There exists $\e_2(n)>0$ with the following property. Assume the hypotheses of Lemma~\ref{lma:Petrunin}, and let $g(t)$ solve \eqref{eqn:RDF} on $M\times(0,T]$ with
\begin{equation}
(1-\e)\wt g(t)\leq g(t)\leq(1+\e)\wt g(t),\qquad 0<\e<\e_2.
\end{equation}
Then there exist $\hat T_2(n,\a)\in(0,1]$ and $C_4(n)>0$ such that
\begin{equation}\label{eqn:single-Morrey}
\int_0^{(T\wedge\hat T_2)r^2}\int_{B_{\wt g(s)}(x_0,r)}
|\wn g|^2\,d\mathrm{vol}_{\wt g(s)}\,ds\leq C_4\e^2r^n.
\end{equation}
\end{lma}
\begin{proof}
It follows from the same argument in \cite[Lemma 4.5]{ChanLaiLee2025} by using Lemma~\ref{lma:Petrunin} instead of the $L^1$ curvature estimate of Petrunin \cite{Petrunin2008}. 
\end{proof}

Similarly, we have the analogous result of \cite[Lemma 4.6]{ChanLaiLee2025}.
\begin{lma}\label{lma:improved-Morrey-grad}
There exists $\e_3(n)>0$ with the following property. Let $\wt g(t)$ be a smooth Ricci flow on $[0,T]$ such that
\begin{enumerate}
\item[(i)] $B_{\wt g(0)}(x_0,8r)\Subset M$;
\item[(ii)] $|\Rm(\wt g(t))|\leq\a t^{-1}$ on $M\times(0,T]$;
\item[(iii)] $\Rm(\wt g(t))+r^{-2}\Id\in\CPIC$ on $M\times[0,T]$.
\end{enumerate}
Suppose $g(t),\hat g(t)$ solve \eqref{eqn:RDF} on $M\times(0,T]$ and
\begin{equation}
\begin{gathered}
(1-\e)\hat g(t)\leq g(t)\leq(1+\e)\hat g(t),\\
(1-\e_3)\wt g(t)\leq g(t),\hat g(t)\leq(1+\e_3)\wt g(t),
\end{gathered}
\end{equation}
with $0<\e<\e_3.$ Then there exist $\hat T_3(n,\a)\in(0,1]$ and $C_5(n)>0$ such that
\begin{equation}\label{eqn:difference-Morrey}
\int_0^{(T\wedge\hat T_3)r^2}\int_{B_{\wt g(s)}(x_0,r)}
|\wn(g-\hat g)|^2\,d\mathrm{vol}_{\wt g(s)}\,ds\leq C_5\e^2r^n.
\end{equation}
\end{lma}
\begin{proof}
The proof follows from the proof of \cite[Lemma 4.6]{ChanLaiLee2025}, using Lemma~\ref{lma:Morrey-grad} and the space-time curvature estimate from Lemma~\ref{lma:Petrunin}.
\end{proof}

\begin{rem}\label{rem:positive-time-background}
For a complete background Ricci flow which is smooth only on $(0,T]$, the preceding integral estimates still hold whenever their curvature assumptions hold for positive time by applying the estimates to the shifted background $\wt g(\tau+t)$, using $|\Rm(\wt g(\tau+t))|\leq\a t^{-1}$, and then letting $\tau\to0$.
\end{rem}

\subsection{Coupling with the heat kernel}

We follow \cite{BCW2019} and use the heat kernel $G(x,t;y,s)$ of
\begin{equation}\label{eqn:G-RF}
\left\{
\begin{aligned}
(\partial_t-\Delta_{\wt g(t)})G(x,t;y,s)&=\mathrm{scal}(\tilde g(t))\, G(x,t;y,s),\\
\lim_{t\downarrow s}G(\cdot,t;y,s)&=\delta_y,
\end{aligned}
\right.
\end{equation}
where the initial delta measure is normalized with respect to $d\mathrm{vol}_{\wt g(s)}$.  In the backward variables it satisfies
\begin{equation}\label{eqn:backward-kernel}
(\partial_s+\Delta_{\wt g(s)})G(x,t;y,s)=0.
\end{equation}

The next Proposition is proved implicitly in \cite[Theorem 4.1]{ChanLaiLee2025}, based on the Li-Yau type estimates of $G(x,t;y,s)$ from  \cite[Proposition 4.8]{ChanLaiLee2025}, see also \cite[Proposition 3.1]{BCW2019}.

Choose $\e_4(n)>0$ sufficiently small that Lemmas~\ref{lma:evo-h-L1} (in its one-solution form), \ref{lma:Morrey-grad}, \ref{lma:improved-Morrey-grad}, and \cite[Lemma 4.7]{ChanLaiLee2025}
 apply whenever the corresponding metrics lie within $\e_4$ of the background. 

\begin{prop}\label{prop:coupled-estimate}
Suppose $\wt g(t)$ is a smooth complete Ricci flow on $M\times(0,T]$ satisfying the assumptions in Theorem~\ref{thm:intro-main}. Let $g(t),\hat g(t)$ solve \eqref{eqn:RDF} on $M\times(0,T]$ and satisfy
\begin{equation}
\begin{gathered}
(1-\e)\hat g(t)\leq g(t)\leq(1+\e)\hat g(t),\\
(1-\e_4)\wt g(t)\leq g(t),\hat g(t)\leq(1+\e_4)\wt g(t),
\end{gathered}
\end{equation}
with $0<\e<\e_4.$ Then there exist $\wt T(n,\a)\in(0,1]$ and $\Lambda_0(n,\a)>0$ such that, for $h=g-\hat g$ and $0<t\leq T\wedge\wt T$,
\begin{equation}\label{eqn:coupled-estimates}
\begin{aligned}
\int_0^t\int_MG(x,t;y,s)|\Rm(\wt g)|(y,s)\,d\mathrm{vol}_{\wt g(s)}(y)\,ds
&\leq\Lambda_0,\\
\int_0^t\int_MG(x,t;y,s)|\wn h|^2(y,s)\,d\mathrm{vol}_{\wt g(s)}(y)\,ds
&\leq\Lambda_0\e^2,\\
\int_0^t\int_M|\wn_yG(x,t;y,s)||\wn h|(y,s)\,d\mathrm{vol}_{\wt g(s)}(y)\,ds
&\leq\Lambda_0\e.
\end{aligned}
\end{equation}
\end{prop}
\begin{proof}
All constants in this proof depend only on $n$ and $\a$. We fix
\begin{equation}
\hat T=\min\{\hat T_1,\hat T_2,\hat T_3,1\},
\end{equation}
and split each time integral into $(0,t/2)$ and $(t/2,t)$.

On $(t/2,t)$, the curvature bound and \cite[Lemma 4.7]{ChanLaiLee2025} give
\begin{equation}
|\Rm(\wt g(s))|\leq\frac{C}{t},\qquad
|\wn h|\leq\frac{C\e}{\sqrt t}.
\end{equation}
Together with \cite[Proposition 4.8]{ChanLaiLee2025}, these prove all three estimates on the later time interval. 

We now explain the covering argument for the earlier time interval with centers fixed in time. The distance distortion estimates under the curvature and Ricci bounds imply
\begin{equation}\label{eqn:distance-distortion}
d_{\tilde g(s)}(p,q)\leq d_{\tilde g(t)}(p,q)+C\sqrt t,\qquad
d_{\tilde g(s)}(p,q)\geq e^{-(n-1)(t-s)}d_{\tilde g(t)}(p,q)
\end{equation}
for $0<s<t$ and $p,q\in M$, see \cite[Lemma 3.1]{SimonToppingGT}. Choose $\b\geq1$ so large that, for $r_0=\b\sqrt t$, we get
\begin{equation}\label{eqn:fixed-to-moving}
B_{\wt g(t)}(z,r_0)\subset B_{\wt g(s)}(z,2r_0),\qquad
\frac t2\leq\hat T(2r_0)^2.
\end{equation}
We restrict $t\leq\wt T$ so that $2r_0\leq1$. At the scale $2r_0$ we have
\begin{equation}
\Rm(\wt g)+(2r_0)^{-2}\Id
=(\Rm(\wt g)+\Id)+\big((2r_0)^{-2}-1\big)\Id\in\CPIC.
\end{equation}

Lemmas~\ref{lma:Petrunin} and \ref{lma:improved-Morrey-grad}, together with Remark~\ref{rem:positive-time-background}, now imply, for every fixed $z\in M$,
\begin{equation}\label{eqn:fixed-ball-L1}
\begin{split}
\int_0^{t/2}\int_{B_{\wt g(t)}(z,r_0)}|\Rm(\wt g(s))|\,d\mathrm{vol}_{\wt g(s)}\,ds
&\leq Ct^{n/2},\\
\int_0^{t/2}\int_{B_{\wt g(t)}(z,r_0)}|\wn h|^2\,d\mathrm{vol}_{\wt g(s)}\,ds
&\leq C\e^2t^{n/2}.
\end{split}
\end{equation}
Volume comparison and \eqref{eqn:fixed-to-moving} give
\begin{equation}\label{eqn:fixed-ball-volume}
\int_0^{t/2}\mathrm{Vol}_{\wt g(s)}\big(B_{\wt g(t)}(z,r_0)\big)\,ds\leq Ct^{n/2+1}.
\end{equation}
By Cauchy-Schwarz, we also obtain
\begin{equation}\label{eqn:fixed-ball-gradient-L1}
\int_0^{t/2}\int_{B_{\wt g(t)}(z,r_0)}|\wn h|\,d\mathrm{vol}_{\wt g(s)}\,ds
\leq C\e t^{(n+1)/2}.
\end{equation}

Consider now the fixed annuli
\begin{equation}
A_k=B_{\wt g(t)}(x,(k+1)r_0)\setminus B_{\wt g(t)}(x,kr_0),
\qquad k=0,1,\ldots.
\end{equation}
A maximal separated set in each annulus gives a cover by $N_k$ balls $B_{\wt g(t)}(z_{kj},r_0)$, with
\begin{equation}\label{eqn:covering-number}
N_k\leq C(k+1)^n\exp\big(C(k+1)r_0\big),
\end{equation}
by volume comparison.  For $y\in A_k$ and $s\leq t/2$, we might assume $\tilde T$ is small so that \cite[Proposition 4.8]{ChanLaiLee2025} and \eqref{eqn:distance-distortion} imply
\begin{equation}
G(x,t;y,s)\leq Ct^{-n/2}e^{-k^2/C},\qquad
|\wn_yG(x,t;y,s)|\leq Ct^{-(n+1)/2}e^{-k^2/C}
\end{equation}
for $y\in A_k$. Applying the first estimate in \eqref{eqn:fixed-ball-L1} on the covering balls and summing over $k$ yields
\begin{equation}
\begin{split}
\int_0^{t/2}\int_MG|\Rm(\wt g)|\,d\mathrm{vol}_{\wt g(s)}\,ds
&\leq C\sum_{k=0}^\infty N_ke^{-k^2/C}\leq C,
\end{split}
\end{equation}
since $r_0\leq1/2$. The second estimate follows in analogously, with an additional factor $\e^2$. For the third, we can use \eqref{eqn:fixed-ball-gradient-L1} and the gradient bound for $G$ to obtain $C\e\sum_kN_ke^{-k^2/C}\leq C\e$. This completes the proof.
\end{proof}

\section{Weak stability of almost WPIC1 Ricci flow}\label{sec:weak-stability}

We now apply the estimates in the previous section to prove our main theorem. We first specify the meaning of the initial data. An $L^\infty$ metric relative to a smooth metric $\wt g_0$ is a measurable positive-definite symmetric two-tensor, defined almost everywhere, which is uniformly equivalent to $\wt g_0$.  
We say that $g(t)$ has {\sl initial data} $g_0$ if
\begin{equation}\label{eqn:initial-trace}
g(t)\longrightarrow g_0\quad\text{in }L^p_{\loc}(M,\wt g_0)
\quad\text{as }t\downarrow0,
\end{equation}
for every $1\leq p<\infty$. In case $g(t)$ is uniformly equivalent to $\tilde g(t)$, it suffices to require local $L^{p_0}$ convergence for some $p_0\in [1,+\infty)$. We restate the statement of our main theorem for the convenience of reader.

\begin{thm}\label{thm:weak-stable-RDF}
For every $n\geq3$ and $\a>0$, there exist $\bar\e_0(n,\a),\tilde T(n,\a)>0$ and $\Lambda(n,\a)\geq1$ with the following property. Let $\wt g(t)$ be a smooth complete Ricci flow on $M\times[0,T]$, with $\wt g(0)=\wt g_0$, satisfying
\begin{enumerate}
\item[(i)] $|\Rm(\wt g(t))|\leq\a t^{-1}$;
\item[(ii)] $\inj(\wt g(t))\geq\sqrt{\a^{-1}t}$;
\item[(iii)] $\Rm(\wt g(t))+\Id\in\CPIC$
\end{enumerate}
on $M\times(0,T]$. If $g_0$ is an $L^\infty$ metric on $M$ such that
\begin{equation}
\|g_0-\wt g_0\|_{L^\infty(M,\wt g_0)}<\bar\e_0,
\end{equation}
then there is a smooth complete solution $g(t)$ to \eqref{eqn:RDF} on $M\times(0,T\wedge \tilde T]$ with initial data $g_0$ and
\begin{enumerate}
\item[(a)]
$\displaystyle \|g(t)-\wt g(t)\|_{L^\infty(M,\wt g(t))}
\leq\Lambda\|g_0-\wt g_0\|_{L^\infty(M,\wt g_0)}$;
\item[(b)] if $g_0$ is continuous on an open set $\Omega\subset M$, then $g(t)\to g_0$ in $C^0_{\loc}(\Omega)$ as $t\downarrow0$.
\end{enumerate}
Moreover, if $g(t),\hat g(t)$ are two solutions with initial data $g_0,\hat g_0$ and
\begin{equation}\label{eqn:stability-class}
\sup_{0<t\leq T\wedge \tilde T}\max\big\{
\|g(t)-\wt g(t)\|_{L^\infty},
\|\hat g(t)-\wt g(t)\|_{L^\infty}\big\}\leq\Lambda\bar\e_0,
\end{equation}
then
\begin{equation}\label{eqn:weak-stability}
\sup_{0<t\leq T\wedge \tilde T}\|g(t)-\hat g(t)\|_{L^\infty(M,\wt g(t))}
\leq\Lambda\|g_0-\hat g_0\|_{L^\infty(M,\wt g_0)}.
\end{equation}
In particular, the solution is unique in this class with its prescribed initial data.
\end{thm}

\begin{proof} 
With Proposition~\ref{prop:coupled-estimate}, the proof of \cite[Theorem 4.1]{ChanLaiLee2025} can now be carried over without  many changes. We sketch the main ideas, focusing on where the proofs differ. Unless otherwise indicated, constants in this proof depend only on $n$ and $\a$.

We first construct a short-time solution starting from $g_0$. We might assume that we are working on the non-compact case. Let $\phi_i$ be a smooth cut-off function on $M$ such that $\phi_i\to 1$. On the support of $\phi_i$, we find a smooth approximation $G_i$ of $g_0$ such that $G_i$ is $\bar\e_0$ close to $g_0$ and  $G_i\to g_0$ in $L^2_{loc}(M,\tilde g_0)$ as $i\to+\infty$. We also let $t_i\to 0$ and consider the smooth metric 
\begin{align*}
    g_{i,0}=\phi_i G_i+(1-\phi_i)\tilde g(t_i)
\end{align*}
on $M$. By the work of Simon \cite{Simon2002}, we might consider the short-time solution $g_i(t)$ on $M\times [0,T_i],$ giving a Ricci-DeTurck flow \eqref{eqn:RDF} with respect to $\tilde g_i(t):=\tilde g(t_i+t)$ and that is $\delta:=\Lambda\bar\e_0$ close to $\tilde g_i(t)$, where $\Lambda(n,\a)$ is a large constant to be chosen. We note that $\tilde g_i(t)$ still satisfies the same assumptions as $\tilde g(t)$. 

Let $h_i=g_i-\tilde g_i$ and suppose that the solution is $\delta$-close to $\wt g$, with $\delta$ small enough for Lemma~\ref{lma:evo-h-L1} and Proposition~\ref{prop:coupled-estimate}. We omit $i$ for convenience.  By Lemma~\ref{lma:evo-h-L1} and \eqref{eqn:backward-kernel}, an integration by parts in the divergence terms gives, after letting $\sigma\downarrow0$ for $v=\sqrt{|h|^2+\sigma}$, 
\begin{equation}\label{eqn:representation}
\begin{split}
|h|(x,t)
&\leq C_1\int_MG(x,t;y,s)|h|(y,s)\,d\mathrm{vol}_{\wt g(s)}(y)\\
&\quad+C_1\int_s^t\int_M G\big[|\wt\Rm||h|^2
+|h||\wn g|^2+|\wn h|(|\wn g|+|\wn\hat g|)\big]\,d\mathrm{vol}_{\wt g(u)}\,du\\
&\quad+C_1\int_s^t\int_M|\wn_yG|
\big[\delta|\wn h|+|h||\wn\hat g|\big]\,d\mathrm{vol}_{\wt g(u)}\,du,
\end{split}
\end{equation}
where $G=G(x,t;\cdot,\cdot )$. We note that the Stokes Theorem in the non-compact case can be justified by the Gaussian bound, see \cite[Lemma 4.9]{ChanLaiLee2025}. 

For $0<\tau\leq T\wedge\wt T$, if we write
\begin{equation}
H(\tau)=\sup_{M\times[0,\tau]}|h|,
\end{equation}
then Proposition~\ref{prop:coupled-estimate} yields
\begin{equation}\label{eqn:absorption}
\begin{split}
H(\tau)
&\leq C_1\sup_{\substack{x\in M\\0\leq t\leq\tau}}\int_MG(x,t;y,0)|h|(y,0)\,d\mathrm{vol}_{\wt g(0)}(y)+C_2\delta H(\tau).
\end{split}
\end{equation}

We choose $\delta$ so that  $C_2\delta\leq1/2$ and thus
\begin{equation}
\sup_{M\times [0,\tau]}|h_i|\leq 2C_1\sup_{\substack{x\in M\\0\leq t\leq\tau}} \int_MG(x,t;y,0)|h_i|(y,0)\,d\mathrm{vol}_{\wt g(0)}(y)\leq C_3\bar\e_0
\end{equation}
using \cite[Proposition 4.8]{ChanLaiLee2025}. If we choose $\Lambda=2C_3$, then it implies $T_i>(T-t_i)\wedge \tilde T$ as $i\to+\infty$. This in particular gives the existence of $g(t)$ by higher order regularity after we let $i\to+\infty$, for example by \cite[Lemma 4.7]{ChanLaiLee2025}. From the analysis above, it also proves (a) and (b) following from a similar argument of \cite[Theorem 5.2]{Simon2002}. 

We now claim that $g(t)$ has initial data $g_0$ in the sense of \eqref{eqn:initial-trace}. Fix $x_0\in M$ and consider $B_{\tilde g_0}(x_0,2)$. For $\e>0$, we choose a smooth metric on $M$ such that $||\tilde g_0-\bar g||_{L^\infty}<\bar\e_0$ and 
\begin{equation}
\int_{B_{\tilde g_0}(x_0,2)}|g_0-\bar g|^2\,d\mathrm{vol}_{\tilde g_0}<\e.
\end{equation}
In the following, we use $L_i$ to denote constants depending also on $\bar g,x_0$ and $\tilde g$. Since 
\begin{equation}
|g-\bar g|\leq |g-\tilde g|+|\tilde g-\tilde g_0|+|\tilde g_0-\bar g|,
\end{equation}
we might assume $|g-\bar g|\leq \delta$ on $B_{\tilde g_0}(x_0,2)\times [0,\hat T]$ for some $\hat T$, by the smoothness of $\tilde g(t)$ up to $t=0$. We denote $\b:=g-\bar g$ so that \eqref{eqn:RDF} implies 
\begin{equation}
\begin{split}
\left(\partial_t -g^{ij}\tilde\nabla_i \tilde\nabla_j \right)\b=g^{-1}*g*\tilde g*\wt\Rm+g^{-1}*g^{-1}*\tilde\nabla g*\tilde\nabla g+g^{-1}*\tilde\nabla^2 \bar g,
\end{split}
\end{equation}
and thus
\begin{equation}
\begin{split}
\left(\partial_t -g^{ij}\tilde\nabla_i \tilde\nabla_j \right) |\b|^2&\leq -|\tilde\nabla \b|^2+C_n|\b| |\tilde\nabla g|^2+L_1\\
&\leq -|\tilde\nabla \b|^2+C_n|\b| \left( |\tilde\nabla \b|^2+|\tilde\nabla \bar g|^2\right)+L_1\\
&\leq -\frac12 |\tilde\nabla \b|^2+L_2
\end{split}
\end{equation}
on $B_{\tilde g_0}(x_0,2)\times [0,\hat T]$. Let $\varphi$ be a time-independent smooth cut-off function such that $\varphi\geq 1$ on $B_{\tilde g_0}(x_0,1)$, vanishing outside $B_{\tilde g_0}(x_0,2)$ and satisfying
$
 |\tilde\nabla \varphi|^2 \leq 10\varphi$. Using again $g=\b+\bar g$ and $|\b|\leq\delta$,
\begin{equation}
\begin{split}
\partial_t \int_M |\b|^2\varphi\,d\mathrm{vol}_{\tilde g(t)}&=\int_M \left(\partial_t -g^{ij}\tilde\nabla_i \tilde\nabla_j \right)|\b|^2\cdot \varphi- \tilde\nabla_i g^{ij}\cdot \tilde\nabla_j |\b|^2\cdot \varphi\,d\mathrm{vol}_{\tilde g(t)}\\
&\quad -\int_M \tilde\nabla_j |\b|^2\cdot \tilde\nabla_i \varphi+\mathrm{scal}(\tilde g) |\b|^2\varphi^2 \,d\mathrm{vol}_{\tilde g(t)}\\
&\leq  \int_M \left(-\frac12 |\tilde\nabla \b|^2+L_2 + C_n|\b||\tilde\nabla g||\tilde\nabla \b|\right) \varphi\,d\mathrm{vol}_{\tilde g(t)}\\
&\quad +\int_ML_2  |\b|^2\varphi+ |\b||\tilde\nabla\b||\tilde\nabla \varphi|\,d\mathrm{vol}_{\tilde g(t)}\\
&\leq -\frac14 \int_M |\tilde\nabla\b|^2 \varphi\,d\mathrm{vol}_{\tilde g(t)}+L_3.
\end{split}
\end{equation}
Hence, 
\begin{equation}
\lim_{t\to 0^+}\int_{B_{\tilde g_0}(x_0,1)}|g(t)-\bar g|^2\,d\mathrm{vol}_{\tilde g(t)}\leq \int_{B_{\tilde g_0}(x_0,2)}|g_0-\bar g|^2\,d\mathrm{vol}_{\tilde g(t)} \leq \e.
\end{equation}
Since $\bar g$ can be made arbitrarily close to $g_0$ in $L^2_{loc}$, the proves the $L^2_{loc}$ convergence as $t\to 0^+$.

The situation with two solutions is similar. Denote $h:g-\hat g$ where both $g$ and $\hat g$ are $\delta$-close to $\tilde g$. Thanks to the $L^2_{loc}$ convergence,  we use the kernel representation as in \eqref{eqn:representation} and Proposition~\ref{prop:coupled-estimate} to get 
\begin{equation}\label{eqn:initial-kernel-control}
\begin{split}
H(\tau)
&\leq C_1\sup_{\substack{x\in M\\0< t\leq\tau}}\limsup_{s\to 0}\int_MG(x,t;y,s)|h|(y,s)\,d\mathrm{vol}_{\wt g(s)}(y)+C_2(\delta+\gamma) H(\tau).
\end{split} 
\end{equation}
where $H(\tau):=\sup_{M\times (0,\tau]} |h|$ now. This is analogous to \eqref{eqn:absorption}. We will fix $\gamma$ so that $C_2\gamma<\frac14$. The same argument above  prove \eqref{eqn:weak-stability}.  

\end{proof}

\begin{rem}\label{rem:trace-needed}
The initial data requirement is part of the uniqueness class. For a measurable metric, convergence on those open sets where it happens to be continuous does not, by itself, specify the initial data. The local $L^p$ formulation above removes this ambiguity without imposing global uniform convergence at time zero.
\end{rem}

\subsection{Initial metrics with a singular set}

We next extend the uniqueness statement to a background which need not be smooth at time zero. 

\begin{thm}\label{thm:unique-general}
Let $M^n$ be a smooth manifold and $(M,\wt d_0,x_0)$ a pointed complete metric space whose metric induces the manifold topology. Suppose $\wt g(t)$ is a smooth complete Ricci flow on $M\times(0,S]$ such that
\begin{enumerate}
\item[(i)] $|{\Rm(\wt g(t))}|\leq\a t^{-1}$;
\item[(ii)] $\inj(\wt g(t))\geq\sqrt{\a^{-1}t}$;
\item[(iii)] $\Rm(\wt g(t))+\Id\in\CPIC$;
\item[(iv)] $d_{\wt g(t)}\to\wt d_0$ locally uniformly on $M\times M$ as $t\downarrow0$;
\item[(v)] there is a closed subset $\mathcal S\subset M$ such that $M\setminus\mathcal S$ is dense, $\mathcal H^n_{\wt d_0}(\mathcal S)=0$, and $\wt g(t)\to\wt g_0$ in $C^0_{\loc}(M\setminus\mathcal S)$ for a continuous Riemannian metric $\wt g_0$ there.
\end{enumerate}
The constants $\Lambda,\bar\e_0$ in Theorem~\ref{thm:weak-stable-RDF} can be chosen so that the following holds. Suppose $g_i(t)$, $i=1,2$, are smooth solutions to \eqref{eqn:RDF} on $M\times(0,S]$ with
\begin{enumerate}
\item[(a)] $\displaystyle\sup_{M\times(0,S]}|g_i-\wt g|\leq\Lambda\bar\e_0$;
\item[(b)] $g_i(t)\to g_{i,0}$ in $C^0_{\loc}(M\setminus\mathcal S)$ as $t\downarrow0$, for continuous metrics $g_{i,0}$ on $M\setminus\mathcal S$.
\end{enumerate}
Then
\begin{equation}\label{eqn:singular-stability}
\sup_{0<t\leq S\wedge \tilde T}\|g_1(t)-g_2(t)\|_{L^\infty(M,\wt g(t))}
\leq\Lambda\|g_{1,0}-g_{2,0}\|_{L^\infty(M\setminus\mathcal S,\wt g_0)}.
\end{equation}
In particular, if $g_{1,0}=g_{2,0}$, then $g_1(t)=g_2(t)$ for every $t\in(0,S]$.
\end{thm}
\begin{proof}
Let us write
\begin{equation}
h=g_1-g_2 \quad \text{ and  } \quad  a_0=\|g_{1,0}-g_{2,0}\|_{L^\infty(M\setminus\mathcal S,\wt g_0)}.
\end{equation}
By Remark~\ref{rem:positive-time-background}, the estimates leading to \eqref{eqn:initial-kernel-control} apply to the set up in the statement above. To get \eqref{eqn:singular-stability},
it remains to control the initial kernel term (i.e. R.H.S. of \eqref{eqn:initial-kernel-control}) without assuming smoothness of the initial data on $\mathcal S$.

By the distance distortion estimate \cite[Lemma 3.1]{SimonToppingGT}, the Ricci lower bound gives
\begin{equation}\label{eqn:distance-to-zero}
d_{\tilde g(s)}(p,q)\leq e^{(n-1)s}\wt d_0(p,q)
\end{equation}
for all $s\in (0,S]$ and $p,q\in M$.
Let $\mu_0$ denote the normalized $n$-dimensional Hausdorff measure of $\wt d_0$, with the normalization for which the Hausdorff measure of a Riemannian metric equals its volume measure. The identity map from $(M,\wt d_0)$ to $(M,d_s)$ is $e^{(n-1)s}$-Lipschitz. Therefore, for every Borel set $E$,
\begin{equation}\label{eqn:measure-domination}
\mathrm{Vol}_{\wt g(s)}(E)%\mu_s(E)
\leq e^{n(n-1)s}\mu_0(E).
\end{equation}
In particular, $\mathrm{Vol}_{\wt g(s)}%\mu_s
(\mathcal S)=0$.

We also need that $\mu_0$ is locally finite. On $M\setminus\mathcal S$, local uniform convergence of the metrics identifies $\mu_0$ with $d\mathrm{vol}_{\wt g_0}$. For a compact subset $D$ of this regular set contained in $B_{\wt d_0}(x_0,R)$, we have
\begin{equation}
\mu_0(D)=\lim_{s\downarrow0}\mathrm{Vol}_{\wt g(s)}%\mu_s
(D)
\leq\limsup_{s\downarrow0}\mathrm{Vol}_{\wt g(s)}%\mu_s
\big(B_{\wt g(s)}(x_0,e^{(n-1)s}R)\big)
\leq C(n,R)
\end{equation}
by \eqref{eqn:distance-to-zero} and volume comparison. Exhausting the regular part of a bounded set by such compact subsets, and using $\mu_0(\mathcal S)=0$, proves the local finiteness of $\mu_0$. 
%This argument also avoids any need to choose a limiting volume density on the singular set.

We fix $x\in M$ and $0<t\leq S\wedge\wt T$. On $M\setminus\mathcal S$, we have by our assumption that
\begin{equation}
|h|(y,s)\longrightarrow|g_{1,0}-g_{2,0}|_{\wt g_0}(y)\leq a_0\, \text{  as } s\to 0^+.
\end{equation}
If $D\subset M$ is compact, \cite[Proposition 4.8]{ChanLaiLee2025}, \eqref{eqn:measure-domination}, and dominated convergence imply
\begin{equation}\label{eqn:singular-initial-local}
\begin{split}
&\int_DG(x,t;y,s)\big(|h|(y,s)-a_0\big)_+\,d\mathrm{vol}_{\wt g(s)}(y)\\
&\qquad\leq Ct^{-n/2}e^{n(n-1)s}
\int_{D\setminus\mathcal S}\big(|h|(y,s)-a_0\big)_+\,d\mathrm{vol}_{\wt g(0)}(y)
\longrightarrow0,
\end{split}
\end{equation}
where the integrands are uniformly bounded by assumption (a) and the local finiteness of $\mu_0$.

To pass from $D$ to $M$, we consider $D=\overline{B_{\wt g(t)}(x,R)}$. On its complement, the distance distortion \eqref{eqn:distance-distortion} gives $d_{\tilde g(s)}(x,y)\geq e^{-(n-1)t}R$ where $0<s<t$. Thus the Gaussian estimates \cite[Proposition 4.8]{ChanLaiLee2025}, the uniform bound condition (a) for $h$ , and volume comparison show that
\begin{equation}
\sup_{0<s<t/2}\int_{M\setminus D}G(x,t;y,s)|h|(y,s)\,d\mathrm{vol}_{\wt g(s)}(y)
\longrightarrow0\quad\text{as }R\to\infty.
\end{equation}
Together with \eqref{eqn:singular-initial-local} and Gaussian estimate of heat kernel, this yields
\begin{equation}\label{eqn:singular-initial-kernel}
\limsup_{s\downarrow0}\int_MG(x,t;y,s)|h|(y,s)\,d\mathrm{vol}_{\wt g(s)}(y)\leq Ca_0.
\end{equation}
Now we may apply \eqref{eqn:initial-kernel-control} to get \eqref{eqn:singular-stability}. If $a_0=0$, the two solutions agree on $(0,S\wedge \tilde T]$. 
\end{proof}

\begin{rem}
The results in Section~\ref{sec:algebra}-\ref{sec:weak-stability} are stated for Ricci flow with almost WPIC1. The methodology indeed only relies on requiring the reference Ricci flow $\tilde g(t)$ satisfies $\tilde R_{1212}+\tilde R_{1313}\geq -2$ for any orthonormal frame $\{e_i\}_{i=1}^n$. We state the result using $\CPIC$ since the non-negativity of $1$-isotropic curvature is known to be preserved under Ricci flow with bounded curvature.
\end{rem}

\section{Existence of expanding solitons by stability}

We say that $(M,g)$ is an expanding gradient Ricci soliton if there is $f\in C^\infty(M)$ such that $2\Ric+g=2\nabla^2f$. It is said to be asymptotic to a metric cone $C(X)$ in the distance sense if its canonical induced flow $g(t)$ satisfies
\begin{equation}
\left(M,g(t),p \right)\longrightarrow C(X)
\end{equation}
as $t\to 0$, in the pointed Gromov-Hausdorff sense.

\subsection{Existence of expander smoothing}

In this sub-section, we show how the stability in Section~\ref{sec:weak-stability} leads to the existence of expanding Ricci solitons, giving an alternative proof to the main result of Deruelle \cite{Deruelle2016}, asserting that a metric cone $C(X)$ over a smooth link $X$ with $\Rm(C(X))\in \mathrm{C}_{\mathrm{PIC2}}$ admits a smoothing. We start with an elementary lemma relating the curvature of the link to that of the cone.
\begin{lma}\label{alg cone lem}
    Let $(X,h)$ be a smooth Riemannian manifold, $C(X)$ the cone over $(X,h)$ with cone metric $g_c:=dr^2+r^2h$. Then the following are equivalent
    \begin{enumerate}
        \item  $\Rm(g_c)\in \CPIC$;
        \item $\Rm(h)-\Id_h\in \mathrm{C}_{\mathrm{PIC2}}$;
        \item  $\Rm(g_c)\in \mathrm{C}_{\mathrm{PIC2}}$.
    \end{enumerate}
\end{lma}
    \begin{proof}
By scaling, it suffices to show $(1)\implies (2)$ and $(2)\implies (3)$ at $r=1$.

For $(1)\implies (2)$, let $\{e_i\}_{i=1}^4$ be any orthonormal four frame of $TX$ ar $r=1$, $\mu$ and $\lambda\in [0,1]$. This induces the following orthonormal frame in $TC(X)\oplus \R$:
\begin{equation}
\left\{
\begin{array}{ll}
\tilde e_1:=e_1, \tilde e_2:=\mu e_2 +\sqrt{1-\mu^2}\partial_t,\\[2mm]
 \tilde e_3:=e_3,\tilde e_4:=\lambda e_4 +\sqrt{1-\lambda^2}\partial_r.
\end{array}
\right.
\end{equation}
where $\partial_t$ denotes the direction of $\R$. 

As $\Rm(g_c)\in \CPIC$ and vanishes on $\text{span}\{ \partial_r\}\oplus \R$, 
\begin{equation}
\begin{split}
0 &\le R(\tilde e_1, \tilde e_3,\tilde e_1,\tilde e_3)+R(\tilde e_1,\tilde e_4,\tilde e_1,\tilde e_4)+R(\tilde e_2,\tilde e_3,\tilde e_2,\tilde e_3)\\
&\quad +R(\tilde e_2,\tilde e_4,\tilde e_2,\tilde e_4)-2R(\tilde e_1,\tilde e_2,\tilde e_3,\tilde e_4)\\
&=R_c(e_1, e_3,e_1,e_3)+\lambda^2 R_c( e_1,e_4,e_1,e_4)+\mu^2 R_c( e_2, e_3,e_2,e_3)\\
&\quad +\lambda^2\mu^2 R_c(e_2,e_4,e_2,e_4)-2\lambda\mu R_c( e_1,e_2,e_3,e_4)
\end{split}
\end{equation}
where $R_c:=\Rm(g_c)$ and $R:=\Rm(g_c\oplus \R)$.

Since $\Rm(g_c)=\Rm(h)-\Id_h$ when restricted on $TX$ at $r=1$, we see that $\Rm(h)-\Id_h\in \mathrm{C}_{\mathrm{PIC2}}$.

\medskip

Next we prove that $(2)\implies (3)$. Let $\{w_i\}_{i=1}^4$ be any orthonormal four frame of $TC(X)\oplus \R^2$. For each $i$, we can decompose $w_i=v_i+y_i$, where $v_i\in TX$ any $y_i\in \text{  span}\{ \partial_r\}\oplus \R^2$. We complexify $R_h:=\Rm(h)$ and $R_c:=\Rm(g_c)$. We  consider  $\eta=v_1+iv_2$, $\xi=v_3+iv_4 \in T_\mathbb{C}X$. By the equivalent condition of $\Rm(h)-\Id_h\in \mathrm{C}_{\mathrm{PIC2}}$ and Bianchi identity, see \cite[Proposition 7.18 (iii)]{Brendle10}, at $r=1$
\begin{equation}
\begin{split}
0&\leq \left(R_h-\mathrm{Id}_h\right)(\eta,\xi,\bar \eta,\bar \xi)=R_c(\eta,\xi,\bar \eta,\bar \xi)\\
&= R_c(v_1, v_3,v_1,v_3)+R_c( v_1,v_4,v_1,v_4)+R_c( v_2, v_3,v_2,v_3)\\
&\quad +R_c(v_2,v_4,v_2,v_4)-2R_c( v_1,v_2,v_3,v_4)\\
&=R(w_1, w_3, w_1,w_3)+R(w_1,w_4,w_1,w_4)+R(w_2,w_3,w_2,w_3)\\
&\quad +R(w_2,w_4,w_2,w_4)-2R(w_1,w_2,w_3,w_4)
\end{split}
\end{equation}
where $R:=\Rm(g_c\otimes\mathbb{R}^2)$. here we used again the fact that $\Rm(g_c)$ vanishes if one of the inputs belongs to $\text{  span}\{ \partial_r\}\oplus \R^2$. Therefore, $\Rm(g_c)\in \mathrm{C}_{\mathrm{PIC2}}$. This completes the proof.
\end{proof}

\medskip

We now show that $C(X)$ admits a smoothing. Note that, under the conditions being considered, $X \simeq \mathbb{S}^{n-1}$. We restate our result for convenience.

\begin{thm} \label{thm:expand-exist}
   Let $(\mathbb{S}^{n-1},h)$ be a smooth closed manifold with $\Rm(h)-\Id_h\in \mathrm{C}_{\mathrm{PIC2}}$. Then there exists an expanding gradient Ricci soliton $(M,g,f,p)$ with $\Rm(g)\in \mathrm{C}_{\mathrm{PIC2}}$  asymptotic to the cone $C(\mathbb{S}^{n-1})$ over $(\mathbb{S}^{n-1},h)$, in the distance sense. Moreover if in addition $\Rm(h)-\Id_h\in \mathrm{int}\left(\mathrm{C}_{\mathrm{PIC2}}\right)$ somewhere on $\mathbb{S}^{n-1}$, then  $\Rm(g)\in\mathrm{int}(\mathrm{C}_{\mathrm{PIC2}})$.
\end{thm}

\begin{proof} 
By the result of Brendle-Schoen \cite{BS09} and strong maximum principle \cite{BS08}, the normalized Ricci flow starting from $h$ will converge to a round metric on $\mathbb{S}^{n-1}$. After suitable scaling and reparametrization, we may suppose that there is a continuous family (in $C^0$ sense) of metrics $\{h_{0,s}\}_{s\in [0,1]}$ on $\mathbb{S}^{n-1}$ such that the following hold
    \begin{itemize}
        \item $h_{0,0}=cg_{\mathbb{S}^{n-1}}$ for some constant $c\in (0,1)$ and $h_{0,1}=h$;
        \item $\Rm(h_{0,s})-I_{h_{0,s}}$ lies inside the interior of $\mathrm{C}_{\mathrm{PIC2}}$ for $s\in [0,1)$.
    \end{itemize}

Consider the metric cone $C(X_s)$ with  the cone metric $g_{c_s}:=dr^2+r^2h_{0,s}$. 
We use the method of continuity as in \cite{Deruelle2016}, except that we replace the PDE method by our stability result. Let $S$ be the subset of $[0,1]$ such that there exists a smooth manifold $M$, a solution $g(t)$ to the Ricci flow on $M\times (0,+\infty)$ and $p\in M$ such that for some $\a>0$,
    \begin{itemize}
    \item $|\Rm(g(t))|\leq \a t^{-1}$;
    \item $\Rm(g(t))\in \mathrm{C}_{\mathrm{PIC2}}$;
    \item $(M,g(t),p)\to C(X_s)$ as $t\to 0$ in pointed Gromov-Hausdorff sense.
\end{itemize}     

We will show that $S=[0,1]$. Assuming this is the case, we let $(M,g(t),p)$ be the pointed Ricci flow corresponding $1\in S$. By \cite[Lemma 3.1]{SimonToppingGT} and Colding's volume continuity, we have $\mathrm{AVR}(g(t))\geq v_0>0$ and hence \cite[Theorem 6.1]{ChanLeePeachey2024} implies that $(M,g(1),f,p)$ is a gradient expanding Ricci soliton for some $f\in C^\infty(M)$. Therefore in order to prove the existence, it suffices to show that $S$ is non-empty, open and closed. Clearly $0\in S$ since  the rotationally symmetric cone over $(\mathbb{S}^{n-1},cg_{\mathbb{S}^{n-1}})$ can be smoothed out by the Bryant expanding soliton with $\Rm>0$.

Now suppose $s_0\in S$ for some $s_0\in [0,1]$ and let  $(M_1, g_1(t),p_1)$ be the corresponding (pointed) Ricci flow satisfying the above properties, where the pointed Gromov-Hausdorff limit as $t\to 0$ is isometric to $C(X_{s_0})$.
By \cite[Lemma 3.1]{SimonToppingGT}, the distance metric $d_{g_1(t)}$ converges to a distance metric $d_0$ on $M_1$ such that $(M_1,d_0,p_1)$ is locally bi-H\"older to $(M_1,g_1(t),p_1)$. By uniqueness of the pointed Gromov-Hausdorff limit, $(M_1,d_0,p_1)$ is isometric to $C(X_{s_0})$. Therefore we might identify $M_1=C(X_{s_0})$. Moreover by 
\cite[Theorem 1.6]{DSS2022b}, $g_1(t)$ extends smoothly to $g_1(0)$ as $t\to 0$ outside the tips. By smooth convergence, $d_{0}=d_{g_1(0)}$ locally outside tips so that $d_{g_1(0)}=d_{g_{c_{s_0}}}$ locally outside tips. Since both metrics are smooth outside tips, Myers–Steenrod theorem implies that the identity is a smooth local Riemannian isometry and thus $g_1(0)=g_{c_{s_0}}$ on $C(X_{s_0})\setminus \{o_{tips}\}$.

 Let $\e>0$ be small and fixed. For all $s\in [0,1]$ which is sufficiently close to $s_0$, let $g_{2,0}:=dr^2+r^2h_{0,s}$. The metric $g_{2,0}$ also a $L^\infty$ metric on $M_1$ such that 
\begin{equation}
    (1-\e) g_{2,0}\leq  g_{1,0}:=dr^2+r^2h_{0,s_0}\leq (1+\e) g_{2,0}.
\end{equation}
As in the proof of \cite[Theorem 5.4]{ChanLaiLee2025}, for $\sigma\to 0^+$, we consider the $L^\infty$ metric 
\begin{equation}
    g_{2,0,\sigma}:=\phi_\sigma g_1(\sigma)+ (1-\phi_\sigma) g_{2,0}
\end{equation}
where $\phi_\sigma(x):=\phi( d_{g_1(\sigma)}(x,o_{tip})/\sigma^{1/4}),$ and $o_{tip}$ is the tip of the cone $C(\mathbb{S}^{n-1})$, and $\phi:[0,+\infty)\to [0,1]$ is a smooth non-increasing function such that $\phi\equiv 1$ on $[0,\frac32]$, vanishes outside $[0,2]$. 

As $\sigma\to 0$, $g_{2,0,\sigma}\to g_{2,0}$ smoothly outside the tips. We apply our stability result Theorem \ref{thm:weak-stable-RDF} with reference $g_1(\sigma +t)$ to $g_{2,0,\sigma}$ and obtain a Ricci-DeTurck flow $\hat g_{2,\sigma}(t)$ with respect to $g_1(\sigma+t)$ starting from $g_{2,0,\sigma}$ such that 
\begin{equation}
   \left (1-\Lambda\e\right) g_{1}(\sigma+t)\leq \hat g_{2,\sigma}(t)\leq \left(1+\Lambda\e\right) g_{1}(\sigma+t)
\end{equation}
on $M_1\times (0,+\infty)$, for all $\sigma\to0^+$. By passing $\sigma\to 0^+$ subsequently and higher order estimate for Ricci Deturck flow \cite[Lemma 4.7]{ChanLaiLee2025}, we obtain a smooth solution to the Ricci-DeTurck flow $\hat g_2(t)$ with respect to $g_1(t)$ on $M_1\times (0,+\infty)$ satisfying 
\begin{equation}\label{eqn:Lip-0}
   \left (1-\Lambda\e\right) g_{1}(t)\leq \hat g_{2}(t)\leq \left(1+\Lambda\e\right) g_{1}(t)
\end{equation}
and $|\Rm(\hat g_2(t))|\leq 2\a t^{-1}$, provided that $\e$ is small.

Furthermore, $\hat g_{2}(t)$ is smooth up to $t=0$  and equals to $g_{2,0}$ away from the $o_{tip}$ by \cite[Proposition 2.2]{ChuLee2025}.  By considering the Ricci-DeTurck ODE, we obtain a Ricci flow $\tilde g_2(t)$ coming out of cone over $(\mathbb{S}^{n-1},h_{0,s})$ in the pointed Gromov-Hausdorff sense, and smoothly up to $t=0$ away from the $o_{tip}$.  Note that the lower bound of the $2$-isotropic curvature is bound by $-C(n,\a)\varepsilon/t$, we may then argue as in \cite[Claim 5.1]{ChanLaiLee2025}  to see that $\Rm(\hat g_{2}(t))\in \mathrm{C}_{\mathrm{PIC2}}$. This shows that $S$ is open.

It remains to show that $S$ is closed. Suppose there is a sequence of $s_i\in [0,1]\cap S$ with  $s_i\to s_0$. We let $(M_i,g_i(t),p_i)$ be the corresponding pointed Ricci flow corresponding to each $s_i$. Since the pointed Gromov-Hausdorff limit of $(M_i,g_i(t),p_i)$  as $t\to 0$ is $C(X_{s_i})$ and $s_i\to s_0$, by \cite[Lemma 3.1]{SimonToppingGT} and Colding's volume continuity we might assume for some $v_0>0$ that $\mathrm{AVR}(g_i(t))\geq v_0$ for all $i$. By Perelman's point-picking \cite[Lemma 3.4]{Lai2019} (see also \cite{Perelman2002, BCW2019}), there is $\a>0$ such that $|\Rm(g_i(t))|\leq \a t^{-1}$ for all $i$.  By \cite[Theorem 4.7]{CheegerGromovTaylor1982} and Hamilton's compactness \cite{HamiltonCPT}, $(M_i,g_i(t),p_i)$ sub-converges to a smooth pointed Ricci flow $(M_\infty,g_\infty(t),p_\infty)$. It also follows from the distance distortion \cite[Lemma 3.1]{SimonToppingGT} and $C(X_{s_i})\to C(X_{s_0})$ that $(M_\infty,g_\infty(t),p_\infty)\to C(X_{s_0})$ as $t\to 0$ in pointed Gromov-Hausdorff sense. This shows that $S$ is closed and hence proves the existence. 

Finally we claim that if, initially, $\Rm(h)-\mathrm{Id}_h$ lies in the interior of $\mathrm{C}_{\mathrm{PIC2}}$, then the corresponding Ricci flow $g(t)$ satisfies $\Rm(g(t))\in \mathrm{int}(\mathrm{C}_{\mathrm{PIC2}})$ for $t>0$. We first note that $M$ is diffeomorphic to $\mathbb{R}^n$ since it is topologically $C(X_1)$. We fix arbitrary $t>0$. Since $C(X_1)$ is non-flat, $g(t)$ cannot be Einstein and is non-symmetric. Since $C(X_1)$ is smooth outside tips, by \cite[Lemma 5.1]{ChanLee2025a} the curvature of $g(t)$ decays to zero at spatial infinity and hence $g(t)$ is irreducible.

By Berger's holonomy classification, its holonomy group can only be either $SO(n)$ or $U(n/2)$ (if in the latter case $n$ is even).  If the holonomy is 
$SO(n)$, then strong maximum principle argument in \cite{BS08} shows that $\Rm(g(t))\in\mathrm{int}(\mathrm{C}_{\mathrm{PIC2}})$. We will see that $U(n/2)$ is not possible. Suppose on the contrary. Then $ g(t)$ is K\"ahler. Since $\Rm(h)-\Id_{h}$ lies in the interior of $\mathrm{C}_{\mathrm{PIC2}},$ it implies that the kernel of $\Ric(g_{c_s})$ is of dimension $1$ on the cone, given by $\partial_r$. As $ g(t)$ converges locally smooth to $g_{c_s}$ away from tip as $t\to 0$, by K\"ahlerity the two smallest eigenvalues of $\Ric(g(t))$ will converge to $0$, contradicting the fact that dim ker$(\Ric(g_{c_s}))=1$.
\end{proof}

\begin{rem}
In Theorem~\ref{thm:expand-exist}, the link is assumed to be smooth. In this case, the distance asymptotics can be improved to smooth asymptotics by \cite{Deruelle2017} and \cite[Theorem 1.6]{DSS2022b}.
\end{rem}

\subsection{Curvature pinching for expanders with continuous links}
\label{subsec:c0-links}

In this sub-section, we want to discuss  the uniqueness of expanders coming out of cones. We show that if an expander $(M,g)$ coming out of a cone $C(X)$ has $\Rm\in \CPIC$, then $\Rm\in \mathrm{C}_{\mathrm{PIC2}}$ if the link $X$ is Reifenberg and is $\mathrm{PIC2}$ limit space. If the link $X$ is $C^2$, the following was proved by \cite{CaoXie2026}.

\begin{thm}\label{thm:c0-PIC2}
Suppose $(M^n,g,f)$, $n\geq3$, is a complete expanding gradient Ricci
soliton with $\Rm(g)\in \CPIC$ and is asymptotic to a metric cone $C(X)$ in the distance sense.  If $X^{n-1}$ is Reifenberg and is a non-collapsed measured Gromov-Hausdorff limit of compact manifold $(X_i^{n-1},h_i)$ with $\Rm(h_i)-\mathrm{Id}_i\in \mathrm{C}_{\mathrm{PIC2}}$, then  $\Rm(g)\in \mathrm{C}_{\mathrm{PIC2}}$.
\end{thm}
\begin{proof}
By Theorem \ref{thm:expand-exist} and its proof, there exists a sequence of gradient expanders $(M_i,g_i,f_i)$ with WPIC2 curvature smoothly asymptotic to $C(X_i)$ converging to an expander $(M_\infty,g_\infty,f_\infty)$ with $\Rm(g_\infty)\in\mathrm{C}_{\mathrm{PIC2}}$ coming out of $C(X)$ in the distance sense. 
Let $g_\infty(t)$ be the canonical Ricci flow generated by the expander $(M_\infty,g_\infty,f_\infty)$, and $g(t)$ the flow generated by $(M^n,g,f)$ respectively. By \cite[Lemma 3.1]{SimonToppingGT}, we may assume both $g(t)$ and $g_\infty(t)$ live on $C(X)=M\simeq M_\infty=\mathbb{R}^n$. Since $\Rm(g_\infty(t))\in \mathrm{C}_{\mathrm{PIC2}}$, by \cite[Proposition 5.4]{DSS2025} at $x\in M$ where $d_c(x,o_{tips})=1$, we have $\Rm(g(x,t))+t^k\cdot\mathrm{Id}\in \mathrm{C}_{\mathrm{PIC2}}$ for $t\in (0,T_k]$, for any $k\in \mathbb{N}$.  By the distance distortion, i.e. \cite[Lemma 3.1]{SimonToppingGT}, and the self similar structure of $g(t)$, this implies $\Rm(g(1))+r^{-2k}(x)\cdot\mathrm{Id}$ for any $k\in \mathbb{N}$ as $x\to\infty$ where $r(x):=d_{g(1)}(x,o_{tips})$. 

We are now in a position to apply the maximum principle argument of Cao-Xie \cite[Lemma 5.1]{CaoXie2026}, with the modification from \cite{CMZ2025}. 
By the strong maximum principle, we may suppose that the scalar curvature of $(M,g,f)$ is positive everywhere, since otherwise $\Ric\equiv 0$ and hence is flat by  \cite[Lemma A.2]{LeeTopping2025}. We let $\ell$ be the lower bound of $\Rm(g)$, i.e. the largest possible number at that point with $\Rm(g)-\ell\cdot  \mathrm{Id}\in \mathrm{C}_{\mathrm{PIC2}}$. Then wherever $\ell$ is nonpositive, by \cite[Lemma 4.2]{BCW2019} and \cite[Lemma 5.1 and Theorem 5.1]{CaoXie2026} we have
\begin{equation}
    \label{eqn:ellscal}
\Delta_f\left(\frac{\ell}{\mathrm{scal}}\right)+2\left\langle\nabla \ln(\mathrm{scal}),\nabla \left(\frac{\ell}{\mathrm{scal}}\right)\right\rangle\leq 2\ell \cdot \frac{ |\Ric(g)|^2}{\mathrm{scal}^2}
\end{equation}
in the sense of barrier, where $\Delta_f:=\Delta +\langle \nabla f,\cdot \rangle$.

By our a-priori estimate, $\ell\ge -r^{-2k}$ as $r\to+\infty$, for all $k\in \mathbb N$. We also have $|\Ric(g)|\leq \mathrm{scal}$ by $\Rm(g)\in \CPIC$ and then we can use the barrier function in \cite[(2.28)]{CMZ2025}, i.e. outside a fixed large compact set  
\[
\Delta_fQ+2\langle\nabla \ln(\mathrm{scal}),\nabla Q\rangle<0,
\]
where $Q:=\exp\left(\frac{a}{\sqrt{f}}\right)f^{-1}\mathrm{scal}^{-1}>0$ for some large positive constant $a$. We may then apply the maximum principle to $\mathrm{scal}^{-1}\cdot {\ell}+\e Q$ as $\e\to0 $ to see that either $\ell\geq 0$ on $M$ or $\mathrm{scal}^{-1}\cdot {\ell}$ attains a negative interior minimum. By \eqref{eqn:ellscal}, the latter case is impossible and thus $\ell\geq 0$ on $M$.
\end{proof}

\begin{cor}\label{cor:uniq}
Suppose $(M,g,f)$ and $(\tilde M,\tilde g,\tilde f)$ are complete expanding gradient Ricci soliton with $\Rm\in \CPIC$ and is asymptotic to same metric cone $C(X)$ in the distance sense.  If $X^{n-1}$ is Reifenberg and is a non-collapsed measured Gromov-Hausdorff limit of compact manifold $(X_i^{n-1},h_i)$ with $\Rm(h_i)-\mathrm{Id}_i\in \mathrm{C}_{\mathrm{PIC2}}$, then $(M,g,f)$ is isometric to $(\tilde M,\tilde g,\tilde f)$. 
\end{cor}
\begin{proof}
    It follows from Theorem~\ref{thm:c0-PIC2} that both solutions satisfy $\Rm\in\mathrm{C}_{\mathrm{PIC2}}$. The uniqueness follows from applying \cite[Theorem 1.1]{FlowerLeeSchulze} (see also \cite{Deruelle2016, Flower2026}) to their canonical Ricci flow. 
\end{proof}


\begin{thebibliography}{99}

\bibitem{Bamler2016}
Bamler, R. H.,
{\sl A Ricci flow proof of a result by Gromov on lower bounds for scalar curvature},
Math. Res. Lett. {\bf 23} (2016), no. 2, 325--337.

\bibitem{BCW2019}
Bamler, R; Cabezas-Rivas, E; Wilking, B.,{\sl  The Ricci flow under almost non-negative
curvature conditions}. Invent. Math. 217 (2019), no. 1, 95–126.

\bibitem{Brendle09}
Brendle, S.,
{\sl A generalization of Hamilton's differential Harnack inequality
for the Ricci flow},
J. Differential Geom. \textbf{82} (2009), no.~1, 207--227.

\bibitem{Brendle10}
Brendle, S., {\sl Ricci flow and the sphere theorem. Graduate Studies in Mathematics},
111. American Mathematical Society, Providence, RI, 2010. viii+176 pp

\bibitem{BS08}
Brendle, S.; Schoen, R. M.,
{\sl Classification of manifolds with weakly 1/4-pinched curvatures}.
Acta Math. 200 (2008), no. 1, 1–13.

\bibitem{BS09} Brendle, S.; Schoen, R. M., {\sl Manifolds with $1/4$-pinched curvature are space forms}, J. Amer. Math. Soc. 22 (2009), no. 1, 287–307.


\bibitem{BurkhardtGuim2019}
Burkhardt-Guim, P.,
{\sl Pointwise lower scalar curvature bounds for $C^0$ metrics via regularizing Ricci flow},
Geom. Funct. Anal. {\bf 29} (2019), no. 6, 1703--1772.

\bibitem{BurkhardtGuim2026}
Burkhardt-Guim, P.,
{\sl Smoothing $L^\infty$ Riemannian metrics with nonnegative scalar curvature outside of a singular set},
Math. Ann. {\bf 395} (2026), Paper No. 57.

\bibitem{CaiWang2026}
Cai, J.-B.; Wang, B.,
{\sl The Ricci-DeTurck flow on complete manifolds},
arXiv:2603.22834 (2026).

\bibitem{CaoXie2026}
Cao, H.-D.; Xie, J.,
{\sl Curvature pinching of asymptotically conical gradient expanding
Ricci solitons},
J. Differential Geom. 134 (2026), no. 3, 269–305. 

\bibitem{ChanLaiLee2025}
Chan, P.-Y.; Lai, Y.; Lee, M.-C.,
{\sl Higher-dimensional flying wing steady Ricci solitons},
arXiv:2510.23005 (2025).

\bibitem{ChanLee2025}
Chan, P.-Y.; Lee, M.-C.,
{\sl Gap theorem on Riemannian manifolds using Ricci flow},
Adv. Math. {\bf 482} (2025), part C, Paper No. 110625.

\bibitem{ChanLee2025a}
Chan, P.-Y.; Lee, M.-C.,
{\sl On the weakly conical expanding gradient Ricci solitons},
arXiv:2503.12416

\bibitem{ChanLeePeachey2024}
Chan, P.-Y.; Lee, M.-C.; Peachey, T. L.,
{\sl Expanding Ricci solitons coming out of weakly PIC1 metric cones},
arXiv:2404.12755

\bibitem{CMZ2025}
Chan, P.-Y.; Ma, Z.; Zhang, Y.,
{\sl Hamilton-Ivey estimates for gradient Ricci solitons.}
Calc. Var. Partial Differential Equations 64 (2025), no. 9, Paper No. 300, 31 pp






\bibitem{CheegerGromovTaylor1982}
Cheeger, J.; Gromov, M.; Taylor, M.,
{\sl Finite propagation speed, kernel estimates for functions of the
Laplace operator, and the geometry of complete Riemannian manifolds},
J. Differential Geom. {\bf 17} (1982), no. 1, 15--53.

\bibitem{ChuLee2025}
Chu, J.; Lee, M.-C.,
{\sl Ricci-DeTurck flow from rough metrics and applications to scalar curvature problems},
J. Funct. Anal. {\bf 289} (2025), no. 2, Paper No. 110916.

\bibitem{ChuLeeWan}
Chu, J.; Lee, M.-C.; Wan, J.,
{\sl Rigidity of positive mass theorem with fast metric decay},
arXiv:2607.17236

 
 
\bibitem{DengZhu15}
Deng, Y.; Zhu, X., {\sl Complete non-compact gradient Ricci solitons with nonnegative Ricci curvature}. 
Math. Z. 279 (2015), no. 1-2, 211–226.



\bibitem{Der14b}
Deruelle, A.,
\emph{Stability of non compact steady and expanding gradient Ricci solitons},
Calc. Var. Partial Differential Equations \textbf{54} (2015),
no.~2, 2367--2405.

\bibitem{Deruelle2016}
Deruelle, A.,
{\sl Smoothing out positively curved metric cones by Ricci expanders}.
Geom. Funct. Anal. 26 (2016), no. 1, 188--249.

\bibitem{Deruelle2017}Deruelle, A., 
{\sl Asymptotic estimates and compactness of expanding gradient Ricci solitons.} 
Ann. Sc. Norm. Super. Pisa Cl. Sci. (5) 17 (2017), no. 2, 485--530.


\bibitem{DeruelleLamm2017}
A.~Deruelle and T.~Lamm,
, {\sl Weak stability of Ricci expanders with positive curvature operator},
Math. Z. \textbf{286} (2017), no.~3--4, 951--985.

%\bibitem{DeruelleSchulze2022}
%Deruelle, A.; Schulze, F.,\emph{ A relative entropy and a unique continuation result for Ricci expanders}. Comm. Pure Appl. Math. 76 (2023), no. 10, 2613–2692.

\bibitem{DSS2022b}
Deruelle, A.; Schulze, F.; Simon, M.,
{\sl On the regularity of Ricci flows coming out of metric spaces}.
J. Eur. Math. Soc. (JEMS) 24 (2022), no. 7, 2233--2277.

\bibitem{DSS2025}
Deruelle, A.; Schulze, F.; Simon, M.,{\sl  Initial stability estimates for Ricci flow and 3-dimensional Ricci-pinched manifolds}. Duke Math. J. 174 (2025), no. 15, 3317–3376.


\bibitem{DeTurck}DeTurck, D. M., {\sl Deforming metrics in the direction of their Ricci tensors}, J. Differential Geom. 18 (1983), no. 1, 157–162.



\bibitem{Flower2026}
Flower, L.,
{\sl On the uniqueness of Ricci flows from certain Alexandrov spaces},
arXiv:2610.03129v1 (2026).


\bibitem{FlowerLeeSchulze}  Flower, L.;  Lee, M.-C.;  Schulze, F., {\sl Uniqueness of Ricci flows from singular initial data}, arXiv:2611.06814 

%\bibitem{Gromov2014}
%Gromov, M.,
%{\sl Dirac and Plateau billiards in domains with corners},
%Cent. Eur. J. Math. {\bf 12} (2014), no. 8, 1109--1156.

\bibitem{Ham82}
Hamilton, R. S.,
{\sl Three-manifolds with positive Ricci curvature}.
J. Differential Geometry 17 (1982), no. 2, 255–306.

\bibitem{HamiltonCPT}Hamilton, R. S., {\sl A compactness property for solutions of the Ricci flow}. Amer. J. Math. 117 (1995), no. 3, 545–572.


%\bibitem{Katsuda1985}
%Katsuda, A.,
%{\sl Gromov's convergence theorem and its application},
%Nagoya Math. J. {\bf 100} (1985), 11--48.

\bibitem{KochLamm2012}
Koch, H.; Lamm, T.,
{\sl Geometric flows with rough initial data},
Asian J. Math. {\bf 16} (2012), no. 2, 209--235.

\bibitem{Lai2019}
Lai, Y.,
{\sl Ricci flow under local almost non-negative curvature conditions}.
Adv. Math. 343 (2019), 353–392.

\bibitem{LammSimon2023}
Lamm, T.; Simon, M.,
{\sl Ricci flow of $W^{2,2}$-metrics in four dimensions},
Comment. Math. Helv. {\bf 98} (2023), no. 2, 261--364.

\bibitem{Lee2026}
Lee, M.-C.,
{\sl Quantification of scalar curvature under $C^0$ convergence using smoothing},
preprint, arXiv:2604.17759.

\bibitem{LeeTopping2025}
Lee, M.-C.; Topping, P. M.,
{\sl Manifolds with PIC1 pinched curvature},
Geom. Topol. {\bf 29} (2025), no. 9, 4767--4798.

\bibitem{LitzingerSimon2025}
Litzinger, F.; Simon, M.,
{\sl Ricci-DeTurck flow of almost continuous $L^2$-metrics, and metrics with distributional scalar curvature bounded from below},
arXiv:2511.23234 (2025).


\bibitem{MantegazzaMascellaniUraltsev}
Mantegazza, C.; Mascellani, G.; Uraltsev, G.,
{\sl On the distributional Hessian of the distance function},
Pacific J. Math. {\bf 270} (2014), no. 1, 151--166.

\bibitem{MW2012}
Munteanu, O.; Wang, J.{\sl  Analysis of weighted Laplacian and applications to Ricci solitons}
. Comm. Anal. Geom. 20 (2012), no. 1, 55–94.


\bibitem{Perelman2002}
Perelman, G.,
{\sl The entropy formula for the Ricci flow and its geometric applications},
arXiv:math.DG/0211159

\bibitem{Petrunin2008}
Petrunin, A. M.,
{\sl An upper bound for the curvature integral},
Algebra i Analiz {\bf 20} (2008), no. 2, 134--148;
English translation in St. Petersburg Math. J. {\bf 20} (2009), no. 2, 255--265.

%\bibitem{Shi1989}
%Shi, W.-X.,
%{\sl Ricci deformation of the metric on complete noncompact Riemannian manifolds},
%J. Differential Geom. {\bf 30} (1989), no. 2, 303--394.

\bibitem{Simon2002}
Simon, M.,
{\sl Deformation of $C^0$ Riemannian metrics in the direction of their Ricci curvature},
Comm. Anal. Geom. {\bf 10} (2002), no. 5, 1033--1074.

\bibitem{SimonToppingGT}
Simon, M.; Topping, P. M.,
{\sl Local mollification of Riemannian metrics using Ricci flow, and Ricci limit spaces},
Geom. Topol. {\bf 25} (2021), no. 2, 913--948.

\bibitem{SimonToppingJDG}
Simon, M.; Topping, P. M.,
{\sl Local control on the geometry in 3D Ricci flow},
J. Differential Geom. {\bf 122} (2022), no. 3, 467--518.

\bibitem{ToppingPIC1}
Topping, P. M.,
{\sl Ricci flow and PIC1},
Surveys in differential geometry 2022: Essays on geometric flows---celebrating 40 years of Ricci flow,
189--211, Surv. Differ. Geom. {\bf 27}, International Press, Somerville, MA, 2024.

\bibitem{ZhaoZhu22}
Zhao, Z.; Zhu, X.,
{\sl Rigidity of the Bryant Ricci soliton},
arXiv:2212.02889.

\end{thebibliography}
\end{document}